\documentclass{ifacconf}

\usepackage{natbib}            % you should have natbib.sty
\usepackage{graphicx}          % Include this line if your 
\usepackage{amsmath,amssymb,color}
\newcommand{\rme}{\mathrm{e}}
\newcommand{\Ddt}{\frac{d}{dt}}

\newcommand{\dds}{\textstyle\frac{d}{ds}\displaystyle}
\newcommand{\vect}[1]{\boldsymbol{#1}}

\newcommand{\pprime}{{\prime\prime}}
\newcommand{\fracdiff}{D_t}
\newcommand{\fracdiffz}{D_z}
\newcommand{\fracdiffc}{{^C}\!D_t} % Caputo
\newcommand{\fracint}{J_t}
\newcommand{\nullng}{0}

\begin{document}

\begin{frontmatter}

\title{An Algebraic Approach for Identification of Linear Systems with Fractional Derivatives} % Title, preferably not more than 10 words.

% \thanks[footnoteinfo]{Financial support by the Deutsche Forschungsgemeinschaft (Ru~538/5) is gratefully acknowledged.}

\author[First]{Nicole Gehring} 
\author[First]{Joachim Rudolph} 
% \author[Third]{Third C. Author}

\address[First]{Chair of Systems Theory and Control Engineering, \\
   Saarland University, 66123 Saarbr\"ucken, Germany \\
   (e-mail: \{n.gehring, j.rudolph\}@lsr.uni-saarland.de).}
% \address[Second]{Colorado State University, Fort Collins, CO 80523 USA (e-mail: author@lamar. colostate.edu)}
% \address[Third]{Electrical Engineering Department, Seoul National University, Seoul, Korea, (e-mail: author@snu.ac.kr)}

%\begin{keyword}                           % Five to ten keywords,  
fractional order systems, fractional derivatives, parameter identification, system identification, algebraic approaches               % chosen from the IFAC 
%\end{keyword}                             % keyword list or with the 
                                          % help of the Automatica 
                                          % keyword wizard

\begin{abstract}                          % Abstract of not more than 250 words.
Identification of fractional order systems is considered from an algebraic point of view. It allows for a simultaneous estimation of model parameters and fractional (or integer) orders from input and output data. It is exact in that no approximations are required. Using Mikusi\'nski's operational calculus, algebraic manipulations are performed on the operational representation of the system. The unknown parameters and (fractional) orders are calculated solely from convolutions of known signals. A generalized Voigt model describing a viscoelastic material is used to illustrate the approach.
\end{abstract}

\end{frontmatter}

% hom und inhom focus on obtaining equation easy to interpret but free of frac deriv
% 
% Since signal derivatives of any order amplify noise, it has to be avoided. Instead in certain cases, (fractional) integration can be used. %which decreases the impact of white noise
% 
% not on how only algebraic identification equation will be written (no time interpretation)
% 
% identification equation in operational notation

\section{Introduction}

Fractional order models have gained increasing interest over the last years. %Whilst a derivative of integer order is of a local nature in the sense that the derivative of a function at a point depends only on nearby values this is not the case for a non-integer, i.e. fractional, derivative at a point which also depends on information of the function in a larger interval. Therefore, it is better suited to account for the some "memory" in certain materials or mediums. %polymeric 
\cite{torvik84} gave one of the first mathematical justifications for the use of such models for viscoelastic materials. However, fractional models have been utilized for a wide spectrum of physical systems including batteries, magnet-suspension systems, electrical circuits, as well as biological and chemical systems -- to name but a few. Several examples, including fractional systems with distributed parameters, can be found in the books \cite{oldham74frac} and \cite{podlubny99book}.

In addition, researchers from different domains have given experimental evidence for the usefulness of fractional models by identifying their parameters and fractional orders. Some of the recently suggested identification procedures use fractional state variable filters based on a known fractional order (\cite{oust01filter}), frequency response functions (\cite{kim09frf}), or finite element methods for the approximation of fractional derivatives (\cite{gaul02femproc}). An overview on system identification for fractional models can be found in \cite{oust07overview}. Nevertheless, most methods known to the authors rely on some kind of approximation.% of differentiation order.

% ECC 2001 Oustaloup: identification of coefficients via state variable filter where fractional orders are considered to be known

An algebraic approach was used in \cite{FliSira03cocv} to identify parameters in ordinary differential equations and in \cite{RudWoi08corsica} for partial differential equations. The present contribution extends the method to linear fractional models, both, with lumped and distributed parameters. The approach allows for the identification of system parameters and fractional orders. It gives exact relations, in the sense that no approximation is required. Furthermore, no assumptions towards commensurability of fractional orders and system stability have to be made.

The basic idea of the method is to use an operational representation of a fractional model (in the sense of Mikusi\'nski), usually described by an equation relating a (known) input and output, and to eliminate at least the unknown non-integer powers of $s$, corresponding to fractional derivatives, by some algebraic manipulation. Here, the focus lies on obtaining an operational equation, the expressions in which can easily be interpreted as functions of time. Unknown quantities are calculated from relations involving only convolutions of known (input and output) signals.

% Compared to integer order systems, fractional models usually require less parameters.

The paper is structured as follows. First, basic background regarding fractional derivatives and Mikusi\'nski's operational calculus is given in section~\ref{sec:defs}. Then, it is discussed how to obtain an equation suitable for the identification of fractional orders and model parameters in linear fractional models. Homogeneous initial conditions are treated in section~\ref{sec:hom}, inhomogeneous ones\footnote{Only classical definitions of initial conditions are considered here. Recent findings by N.~Maamri and J.C.~Trigeassou as well as T.T.~Hartley and C.F.~Lorenzo (\cite{trigmaam12initial}, \cite{lorhart08initial}) show that these may not be sufficient, in general.} in section~\ref{sec:inhom}. A generalized Voigt model for a viscoelastic material serves as an example. In section~\ref{sec:pde} the results are extended to a fractional distributed parameter system, the diffusion-wave equation of fractional order, followed by a brief conclusion.

% fractional derivations $\dds^\alpha$ are to be avoided, since no correspondence is known for this operation; see Oldham 1974, p.~135 in the context of the Laplace transform $\Rightarrow$ allowed are integer order derivations and multiplications with $s^\alpha$, $\alpha\in\Rset$

% Scott-Blair stress-strain law with fractional derivatives to introduce a material property that is intermediate between the elastic modulus (Hooke solid) and the coefficient of viscosity (Newton fluid) -- see \cite{mainardi97frac}

% dashpot

\section{General definitions and notational aspects}
\label{sec:defs}

In this section the most important basic definitions used in this paper are revisited without being exhaustive on technical details. For a elaborate introduction into the mathematical background of fractional systems the reader is referred to \cite{oldham74frac} and \cite{podlubny99book} or for a well written brief summary to \cite{gorenflo97frac}. The calculus of Mikusi\'nski is detailed in \cite{MikOpCal}.

\subsection{Fractional derivatives}

% With fractional differentiation the order of differentiation is generalized to fractional orders. Several definition are known.
One of the (or maybe the) most commonly used definition(s) for fractional derivatives is the one attributed to B.~Riemann and J.~Liouville. A fractional derivative of order\footnote{Zero is included in the set $\Rset_+$.} $\alpha\in\Rset_+$ of %an integrable
a function $f$ is defined as %in \cite{gorenflo97frac}:
\begin{equation}
\label{eq:fracdiff_def}
  \fracdiff^\alpha f(t) = \begin{cases}
    \frac{1}{\Gamma(\nu-\alpha)}\Ddt^\nu\int_0^t\frac{f(\sigma)}{(t-\sigma)^{\alpha+1-\nu}}\,d\sigma, &\alpha\neq\nu, \\[2ex]
    \displaystyle\Ddt^\nu f(t), &\alpha=\nu,
  \end{cases}
\end{equation}
$t>0$, where $\nu-1<\alpha\le\nu$, $\nu\in\Nset$ and $\Gamma$ denotes the Gamma function. The inverse operation of fractional differentiation is fractional integration. For $t>0$ and $\alpha\in\Rset_+$ it is defined as
\begin{equation}
\label{eq:fracint}
  \fracint^\alpha f(t) = \frac{1}{\Gamma(\alpha)}\int_0^t(t-\sigma)^{\alpha-1}f(\sigma)\,d\sigma.
\end{equation}

An alternative definition has been introduced by M.~Caputo in 1990 and is often referred to as Caputo fractional derivative:
\begin{equation}
\label{eq:fracdiffc_def}
  \fracdiffc^\alpha f(t) = \begin{cases}
    \dfrac{1}{\Gamma(\nu-\alpha)}\displaystyle\int_0^t\dfrac{f^{(\nu)}(\sigma)}{(t-\sigma)^{\alpha+1-\nu}}\,d\sigma, &\alpha\neq\nu \\[2ex]
    \displaystyle\Ddt^\nu f(t), & \alpha=\nu,
  \end{cases}
\end{equation}
$\nu-1<\alpha\le\nu$, $\nu\in\Nset$. In contrast to the Riemann-Liouville definition it is not necessary to define fractional order initial conditions, making Caputo's definition more suitable in the context of solving equations with fractional derivatives. %This fact is illustrated by the following equivalences
% \begin{subequations}
%   \begin{align}
%     \fracdiff^\alpha f(t) &= \fracdiff^\alpha g(t) \quad\Leftrightarrow\quad f(t) = g(t) + \sum_{j=1}^\nu c_j t^{\alpha-j} \\
%     \fracdiffc^\alpha f(t) &= \fracdiffc^\alpha g(t) \quad\Leftrightarrow\quad f(t) = g(t) + \sum_{j=1}^\nu d_j t^{\nu-j},
%   \end{align}
% \end{subequations}
% $\nu-1<\alpha\le\nu$, $\nu\in\Nset$, where $c_j$ and $d_j$ are arbitrary coefficients.

Note that other definitions of fractional derivatives are known, like the one due to A.K.~Gr\"unwald and A.V.~Letnikov which is especially useful when dealing with discrete approximations (e.g.~\cite{gaul02femproc}). For a wide class of functions the Riemann-Liouville and the Gr\"unwald-Letnikov definition are equivalent (see~\cite{podlubny99book}). %Hence, the latter can be used for obtaining a numerical solution whilst the first is more suitable for calculation.

In literature on parameter identification, most researchers seem to use the definition \eqref{eq:fracdiff_def} due to Riemann and Liouville (e.g.~\cite{torvik84,bagley83,oldham74frac,oust07overview}). However, in most of the cases homogeneous initial conditions are assumed where \eqref{eq:fracdiff_def} and \eqref{eq:fracdiffc_def} are equivalent. Here, both definitions will be treated for homogeneous and inhomogeneous initial conditions.

\subsection{Operational calculus}

In this paper Mikusi\'nki's operational calculus is used (see \cite{MikOpCal}). A brief introduction of this calculus in the context of fractional systems can be found in \cite{hotzel98theory}. However, readers unfamiliar with this calculus can (in a simplified manner) consider the operational expressions as Laplace transforms.

% though the calculus is more restrictive in that context (e.g.\ differentiability of $f$ required).

Using Mikusi\'nski's operational calculus \eqref{eq:fracdiff_def} reads
\begin{subequations}
  \begin{align}
  \label{eq:fracdiff_op}
    \fracdiff^\alpha f &= \left[s^\alpha\hat f - \sum_{k=0}^{\nu-1}\fracdiff^k\fracint^{\nu-\alpha}f(\nullng)s^{\nu-1-k}\right] \\
\intertext{and \eqref{eq:fracdiffc_def} yields}
  \label{eq:fracdiffc_op}
    \fracdiffc^\alpha f &= \left[s^\alpha\hat f -\sum_{k=0}^{\nu-1}f^{(k)}(\nullng)s^{\alpha-1-k}\right],
  \end{align}
\end{subequations}
$\nu-1<\alpha\le\nu$, $\nu\in\Nset$, where $\fracint^{\nu-\alpha}f(\nullng)$ is the limit of $\fracint^{\nu-\alpha}f(t)$ for $t\rightarrow\nullng$. Note that both correspondences are equivalent in the case of homogeneous initial conditions. Hence, for \eqref{eq:fracint} it follows
\begin{equation}
  \fracint^\alpha f = \left[s^{-\alpha}\hat f\right], \qquad \alpha\in\Rset_+.
\end{equation}
% laplace transform in \cite[p.~134]{oldham74frac}

Two fundamental properties of this operational calculus are
\begin{equation}
\label{eq:corres_falt}
  [\hat f_1\hat f_2](t)=(f_1\!\star\!f_2)(t)=\int_0^t f_1(\sigma)f_2(t-\sigma)d\sigma
\end{equation}
and
\begin{equation}
\label{eq:corres_derv}
  [\hat f^\prime](t)=-(\bar tf)(t)=-tf(t)
\end{equation}
where $\bar t:t\mapsto t$ is the identity map and $\hat f^\prime=\dds\hat f$ denotes the derivative w.r.t.\ $s$.

\section{The case of homogeneous initial conditions}
\label{sec:hom}

Based on the fundamental notions above, in this section, the identification problem is addressed for fractional systems assuming homogeneous initial conditions. This way, all basic ideas can later quite easily be adapted to the general case.

\subsection{An introductory example}

A commonly used example is the empirical model
\begin{equation}
\label{eq:ex_visco}
  \sigma(t) = E_0\varepsilon(t) + E_1\fracdiff^\alpha\varepsilon(t)
\end{equation}
of a viscoelastic material\footnote{A material is considered elastic for $\alpha=0$ and viscous for $\alpha=1$.}, also referred to as three parameter generalized Voigt model (e.g.~\cite{bagley83,podlubny99book}). The stress $\sigma$ is expressed as a superposition of an elastic part $E_0\varepsilon$ and a viscoelastic part $E_1\fracdiff^\alpha\varepsilon$ using a fractional derivative of order $\alpha\in(0,1)$ of the strain $\varepsilon$. 

The aim here is to identify the parameters $E_0$, $E_1$, as well as the fractional order $\alpha$ from known signals $\sigma$ and $\varepsilon$. Therefore, using definition \eqref{eq:fracdiff_op} for an operational notation, the fractional model \eqref{eq:ex_visco} is written as
\begin{equation}
\label{eq:ex_visco_1}
  \hat\sigma = E_0\hat\varepsilon + E_1s^\alpha\hat\varepsilon,
\end{equation}
where for simplicity a homogeneous initial condition $\varepsilon(0)=0$ is assumed.

In order to obtain an equation that can easily be interpreted and that does not involve any fractional derivatives, $\dds$ is applied to \eqref{eq:ex_visco_1}:
\begin{equation}
\label{eq:ex_visco_2}
  \hat\sigma^\prime = E_0\hat\varepsilon^\prime + E_1(\hat\varepsilon^\prime+\alpha s^{-1}\hat\varepsilon)s^\alpha.
\end{equation}
Then a combination of \eqref{eq:ex_visco_1} and \eqref{eq:ex_visco_2} yields an expression without $s^\alpha$:
\begin{equation}
\label{eq:ex_visco_3}
  (\hat\varepsilon^\prime+\alpha s^{-1}\hat\varepsilon)(\hat\sigma-E_0\hat\varepsilon) = \hat\varepsilon(\hat\sigma^\prime-E_0\hat\varepsilon^\prime).
\end{equation}
Reinterpreting the expressions as functions of time gives
\begin{equation}
\label{eq:ex_visco_4}
  (\varepsilon\!\star\!\bar t\sigma-\bar t\varepsilon\!\star\!\sigma)(t) 
  = -\alpha\!\int_0^t(\varepsilon\!\star\!\sigma)(\tau)\,d\tau + \alpha E_0\!\int_0^t(\varepsilon\!\star\!\varepsilon)(\tau)\,d\tau
\end{equation}
(cf.~\eqref{eq:corres_falt} and \eqref{eq:corres_derv}). Note that neither fractional derivatives nor derivatives of integer order appear in \eqref{eq:ex_visco_4}. In order to calculate $E_0$ and $\alpha$ from trajectories $t\mapsto\varepsilon(t)$ and $t\mapsto\sigma(t)$ at least one more (independent) equation in the parameters is required.

Eq.~\eqref{eq:ex_visco_3} together with the equation obtained by multiplying it with $s^{-1}$ yields a linear system for $\alpha$ and $\alpha E_0$:
\begin{multline}
\label{eq:ex_visco_5}
  \begin{pmatrix} -\int_0^t(\varepsilon\!\star\!\sigma)(\tau)\,d\tau & \int_0^t(\varepsilon\!\star\!\varepsilon)(\tau)\,d\tau \\[.5ex] -\int_0^t\!(t\!-\!\tau)(\varepsilon\!\star\!\sigma)(\tau)\,d\tau & \int_0^t\!(t\!-\!\tau)(\varepsilon\!\star\!\varepsilon)(\tau)\,d\tau \end{pmatrix}\!
  \begin{pmatrix} \alpha \\[.5ex] \alpha E_0 \end{pmatrix} \\
  = \begin{pmatrix} (\varepsilon\!\star\!\bar t\sigma-\bar t\varepsilon\!\star\!\sigma)(t) \\[.5ex] \int_0^t(\varepsilon\!\star\!\bar t\sigma-\bar t\varepsilon\!\star\!\sigma)(\tau)\,d\tau\end{pmatrix},
\end{multline}
where double integrals have been replaced by simple ones using the general relation
\begin{equation}
\label{eq:ex_visco_6}
  \int_0^t\int_0^{\tau_1}\cdots\int_0^{\tau_{k-1}}y(\tau_{k-1})\,d\tau_{k-1} = \int_0^t\frac{(t-\tau)^{k-1}}{(k-1)!}y(\tau)\,d\tau
\end{equation}
for $k$-times integrals, $k>1$. Eq.~\eqref{eq:ex_visco_5} can be solved for the parameters as long as the matrix in the relation is regular. Obviously the matrix is singular for zero signals $\varepsilon$ and $\sigma$. For all other signals the determinant is expected to vanish only at singular points. %However, further assumptions concerning the signals have been made in the previous calculations. When eliminating $s^\alpha$ using \eqref{eq:ex_visco_1} and \eqref{eq:ex_visco_2} both $\hat\varepsilon\neq0$ and $\hat\varepsilon^\prime+\alpha s^{-1}\hat\varepsilon\neq0$ was presumed. The first condition has already been treated and corresponds to $\varepsilon$ not being the zero function. For the second one $\varepsilon(t)\neq Kt^{\alpha-1}$ with arbitrary $K\in\Rset$ is required. Note, that the absolute of the fractional derivative $\fracdiff^\alpha$ of such a function $\varepsilon$ is infinite. \emphng{FALSCH, vgl.\ Determinante im allgemeinen Fall}
  
Based on the parameter estimates obtained from \eqref{eq:ex_visco_5}, the remaining coefficient $E_1$ can easily be calculated by means of \eqref{eq:ex_visco}. However, since differentiations of $\varepsilon$ are to be avoided, \eqref{eq:ex_visco} is integrated (corresponding to a multiplication of \eqref{eq:ex_visco_1} with $s^{-1}$). Solving for the missing parameter then gives
\begin{equation}
\label{eq:ex_visco_7}
  E_1 = \frac{\int_0^t\sigma(\tau)\,d\tau-E_0\int_0^t\varepsilon(\tau)\,d\tau}{\fracint^{1-\alpha}\varepsilon(t)}.
\end{equation}

\begin{figure}[b]
  \vspace{-1ex}
  \includegraphics[width=\linewidth]{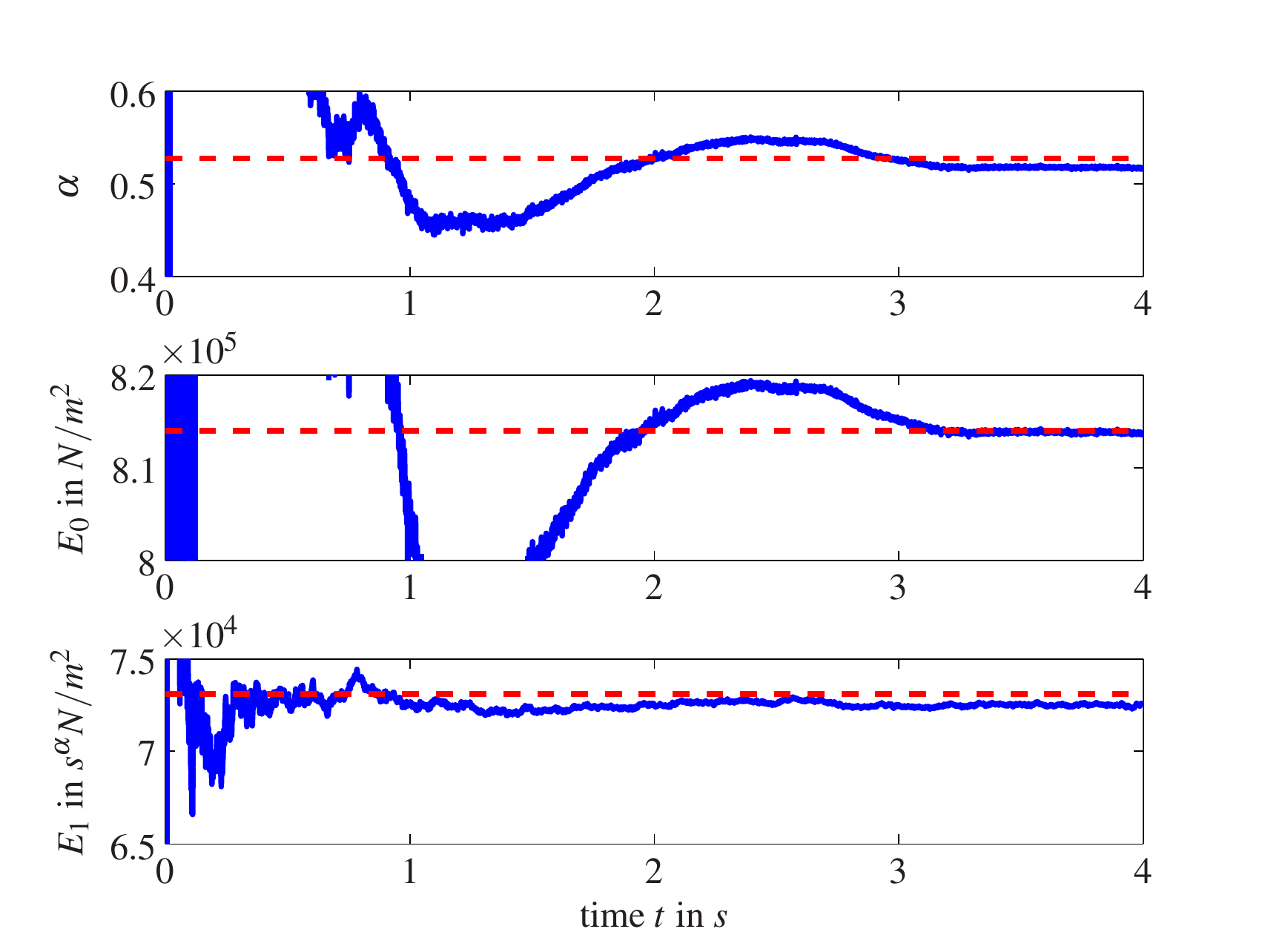}
  \vspace{-3ex}
  \label{fig:ident}
  \caption{Identified values for $\alpha$, $E_0$, and $E_1$ (solid) compared to model parameters (dashed).}
\end{figure}

\begin{figure}[b]
  \vspace{-1ex}
  \includegraphics[width=\linewidth]{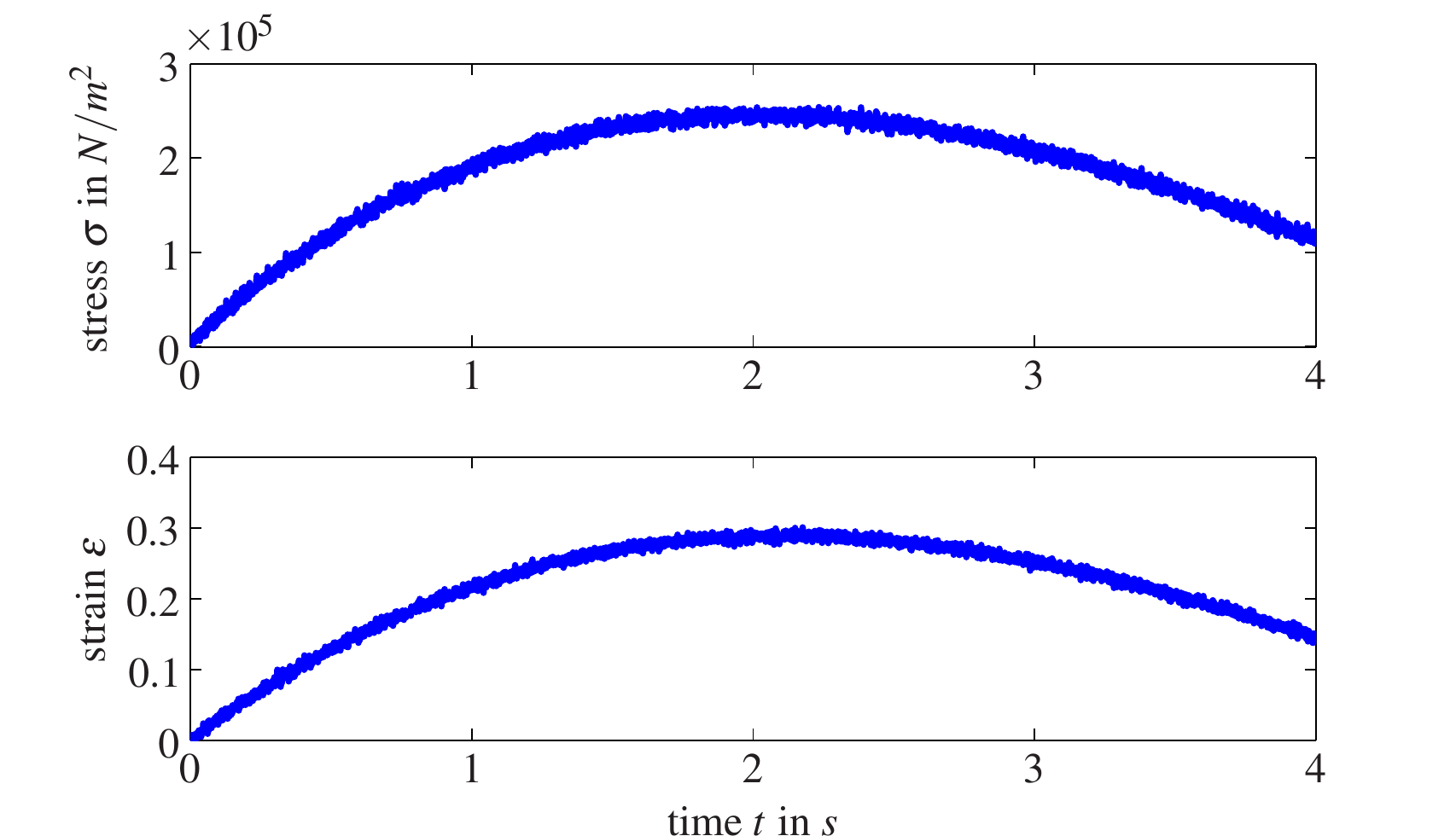}
  \vspace{-3ex}
  \label{fig:signals}
  \caption{Trajectories of stress $\sigma$ and strain $\varepsilon$.}
\end{figure}

For validation purposes a simulation of the introductory example \eqref{eq:ex_visco} was performed with values for Polybutadiene as in \cite{bagley83}. %In the implementation of the system an approximation for fractional derivative arising from the Gr\"unwald-Letnikov definition was used\footnote{Recall that for most function this definition can be considered being equivalent to the one due to Riemann and Liouville.} (see \cite{podlubny99book,gaul02femproc}):
% \begin{equation}
%   \fracdiff^\alpha f(t) \approx \left(\frac{t}{N}\right)^{-\alpha}\sum_{k=0}^{N-1} A_{k+1} f(t-k\tfrac{t}{N})
% \end{equation}
% with $N\in\Nset$ and
% \begin{equation}
%   A_{k+1} = \frac{\Gamma(k-\alpha)}{\Gamma(-\alpha)\Gamma(k+1)} = \frac{k-1-\alpha}{k} A_k.
% \end{equation}
The parameters $E_0$ and $E_1$ as well as the fractional order $\alpha$ are identified using \eqref{eq:ex_visco_5} and \eqref{eq:ex_visco_7}. The integrals in \eqref{eq:ex_visco_7} are approximated using
\begin{equation}
  \fracint^\alpha f(t) \approx \left(\frac{t}{N}\right)^{\alpha}\sum_{k=0}^{N-1} A_{k+1} f(t-k\tfrac{t}{N})
\end{equation}
with $N\in\Nset$ and
\begin{equation}
  A_{k+1} = \frac{\Gamma(k+\alpha)}{\Gamma(\alpha)\Gamma(k+1)} = \frac{k-1+\alpha}{k} A_k,
\end{equation}
which follows from the Gr\"unwald-Letnikov definition of fractional integral and derivatives (cf.~\cite{gaul02femproc}). Fig.~\ref{fig:ident} shows exemplary results of the identification when using $N=4001$ and signals with white noise as depicted in fig.~\ref{fig:signals}. For $t>3$~s, all parameters are identified with an error of less than two percent. The results can be improved even further by integrating \eqref{eq:ex_visco_5} to reduce the impact of the white noise, and by considering additional equations and solving a least squares problem.

\subsection{The general approach}
\label{sec:hom_gen}

In the case of homogeneous initial conditions the definitions \eqref{eq:fracdiff_def} and \eqref{eq:fracdiffc_def} of fractional derivatives are equivalent. Therefore, without loss of generality the one due to Riemann and Liouville is used in this section.

Consider the following relation between two known signals\footnote{The generalization towards multiple equations with potentially more than two signals is obvious.} $u$ and $y$:
\begin{multline}
\label{eq:hom_sys}
  \big(\bar a_0 + \bar a_1\fracdiff^{\bar\alpha_1} + \cdots + \bar a_n\fracdiff^{\bar\alpha_n}\big)y(t) \\
  = \big(\bar b_0 + \bar b_1\fracdiff^{\bar\alpha_{n+1}} + \cdots + \bar b_m\fracdiff^{\bar\alpha_{n+m}}\big)u(t).
\end{multline}
At least some of the parameters $\bar a_i\in\Rset$, $i=0,\dots,n$ and $\bar b_i\in\Rset$, $i=0,\dots,m$ as well as the fractional orders $\bar\alpha_i\in\Rset_+$, $i=1,\dots,n+m$ are considered unknown and are to be identified. Note that if all coefficients in \eqref{eq:hom_sys} are unknown, only the coefficients in a normalized equation can be identified. Define by $\Theta$ the set of unknown parameters.

With \eqref{eq:fracdiff_op}, the operational representation of \eqref{eq:hom_sys} reads
\begin{multline}
\label{eq:hom_sys_1}
  \left(\bar a_0 + \bar a_1s^{\bar\alpha_1} + \cdots + \bar a_ns^{\bar\alpha_n}\right)\hat y \\
  = \left(\bar b_0 + \bar b_1s^{\bar\alpha_{n+1}} + \cdots + \bar b_ms^{\bar\alpha_{n+m}}\right)\hat u.
\end{multline}
First, the fractional operators $s^{\bar\alpha_i}$, $i=1,\dots,m+n$ are eliminated. For that, \eqref{eq:hom_sys_1} is rewritten as
\begin{multline}
\label{eq:hom_sys_2}
  \big(a_0(s)\hat y-b_0(s)\hat u\big)s^{\alpha_0} + \big(a_1(s)\hat y-b_1(s)\hat u\big)s^{\alpha_1} \\
  + \cdots + \big(a_r(s)\hat y-b_r(s)\hat u\big)s^{\alpha_r} = 0
\end{multline}
where $a_i(s),b_i(s)\in\Rset(s)$, $i=0,\dots,r$ with $r\le n+m$ are rational in $s$. These expressions obviously involve the original parameters $\bar a_i$, $i=0,\dots,n$ and $\bar b_i$, $i=0,\dots,m$.  %\emphng{nachvollziehbar?} In this step, knowledge concerning the fractional orders can be used (but is not necessary). Hence $r\le n+m$.
Note that one or more of the fractional orders $\alpha_i$, $i=0,\dots,r$ can be integer or even zero.

As in the introductory example, \eqref{eq:hom_sys_2} is derived w.r.t.\ $s$ a total of $r$ times. The resulting homogeneous linear system of equations for $s^{\alpha_i}$, $i=0,\dots,r$ has the general structure
\begin{multline}
\label{eq:hom_sys_3mat}
  \begin{pmatrix}
    a_0(s)\hat y-b_0(s)\hat u & \cdots & a_r(s)\hat y-b_r(s)\hat u \\[.5ex]
    p_{1,0}\!(s,\!\dds)\hat y\!-\!q_{1,0}\!(s,\!\dds)\hat u & \cdots & p_{1,r}\!(s,\!\dds)\hat y\!-\!q_{1,r}\!(s,\!\dds)\hat u \\[.5ex]
    \vdots & & \vdots \\[.5ex]
    p_{r,0}\!(s,\!\dds)\hat y\!-\!q_{r,0}\!(s,\!\dds)\hat u & \cdots & p_{r,r}\!(s,\!\dds)\hat y\!-\!q_{r,r}\!(s,\!\dds)\hat u
  \end{pmatrix} \\
  \cdot\begin{pmatrix} s^{\alpha_0} \\[.5ex] \vdots \\[.5ex] s^{\alpha_r} \end{pmatrix} = \vect 0
\end{multline}
where $p_{i,j}(s,\dds),q_{i,j}(s,\dds)\in\Rset(s)[\dds]$, $i=1,\dots,r$, $j=0,\dots,r$ are polynomials in $\dds$ with coefficients in $\Rset(s)$. They are defined recursively by
\begin{subequations}
\label{eq:hom_sys_3}
  \begin{equation}
    p_{i,j}(s,\dds) = \dds p_{i-1,j}(s,\dds) + p_{i-1,j}(s,\dds)(\dds+\alpha_js^{-1})
  \end{equation}
  \begin{equation}
    q_{i,j}(s,\dds) = \dds q_{i-1,j}(s,\dds) + q_{i-1,j}(s,\dds)(\dds+\alpha_js^{-1})
  \end{equation}
\end{subequations}
with $p_{0,j}(s,\!\dds)=a_j(s)$ and $q_{0,j}(s,\!\dds)=b_j(s)$, $j=0,\dots,r$. Note that products $p_{i,j}(s,\!\dds)\hat y$ and $q_{i,j}(s,\!\dds)\hat u$ can also be written as sums of $\hat y$ and $\hat u$ and its derivatives with coefficients in $\Rset(s)$. Hence, the operator $\dds$ does not explicitly occur.

From \eqref{eq:hom_sys_3mat} it becomes apparent that the matrix in this equation, henceforth denoted by $P$, has to be singular\footnote{The matrix $P$ is the result of applying the operators, the entries of which belong to a commutative ring.}. Then, from $\det P=0$ an equation is obtained that involves
\begin{itemize}
  \item (most importantly) only integer orders of $s$,
  \item products of $r$ operators $\hat u$, $\hat y$, and their derivatives w.r.t.\ $s$,
  \item as well as products of parameters $\bar a_i\in\Rset$, $i=0,\dots,n$ and $\bar b_i\in\Rset$, $i=0,\dots,m$ and fractional orders $\alpha_i\in\Rset_+$, $i=0,\dots,r$.
\end{itemize}
A multiplication with $s^{-k}$ with sufficiently large $k$ such that only non-positive powers of $s$ remain gives an expression of the form
\begin{equation}
\label{eq:hom_sys_4}
  F(s^{-1},\hat y,\dots,\hat y^{(r)},\hat u,\dots,\hat u^{(r)}) = 0,
\end{equation}
where $F$ is polynomial w.r.t.\ to all its arguments and parameters $\Theta_1\subset\Theta$ including all fractional orders $\alpha_i$, $i=0,\dots,r$ to be identified. The knowledge of at least one of the orders is required in order to identify the others (cf.~\eqref{eq:hom_sys_2}) -- again, this can be considered as a kind of non-restrictive normalization.

Note that the expressions occurring in \eqref{eq:hom_sys_4} can easily be written as functions of time $t$ again. Hence, depending on the number of unknown (independent) coefficients in \eqref{eq:hom_sys_4}, which is larger or equal to the cardinality of $\Theta_1$, additional equations are required to solve for the parameters. They are obtained, e.g., by multiplication of \eqref{eq:hom_sys_4} with negative powers of $s$. Further integration of the corresponding time functions yields independent equations for sufficiently rich signals. While other approaches of generating independent equations are possible -- basically by applying other operators from $\Rset(s)[\dds]$ to \eqref{eq:hom_sys_4} -- which might exhibit a better numerical behavior or might be more suitable in the presence of noise, this is not in the focus of the present contribution. Also, the issue of whether the parameter problem arising from \eqref{eq:hom_sys_4} is solved as a nonlinear equation or a linear one (by overparametrization) is not discussed here. Several approaches are well-known to the community. In the context of algebraic identification some insight into this is given in \cite{FliSira03cocv} and \cite{geh12sysid}.

It is rather obvious that parameters $\Theta_2\subset\Theta$ with $\Theta=\Theta_1\cup\Theta_2$ exclusively occurring as a factor of a fractionally differentiated signal in \eqref{eq:hom_sys} or more precisely a common factor of any $a_i(s)\hat y-b_i(s)\hat u$, $i=0,\dots,r$ in \eqref{eq:hom_sys_2} will not occur in \eqref{eq:hom_sys_4}. However, as done for the introductory example, these parameters can be calculated based on the identified values for $\Theta_1$ and \eqref{eq:hom_sys_2}.

To this end, based on estimates for the fractional orders, \eqref{eq:hom_sys_2} is multiplied with $s^{-\nu}$, $\Nset\ni\nu>\max_i\alpha_i$, i.e., such that only non-positive powers of $s$ remain. Then, as before, depending on the number of parameters remaining to be identified (which corresponds to the cardinality of $\Theta_2$) additional equations are obtained by multiplication with negative powers of $s$. The solution of this linear system of equations completes the identification.

% signal abh\"angigkeit ay+by

\section{The case of inhomogeneous initial conditions}
\label{sec:inhom}

For inhomogeneous initial conditions the definitions \eqref{eq:fracdiff_def} and \eqref{eq:fracdiffc_def} are essentially different. To illustrate this, first, the introductory example is revisited. %One basically has two options: The unknown initial conditions can be considered as additional parameters to be identified or they can be eliminated by algebraic manipulations.

\subsection{Introductory example}

\subsubsection{Riemann-Liouville fractional derivative}

For inhomogeneous initial conditions the operational representation of \eqref{eq:ex_visco} (with the Riemann-Liouville operator) reads
\begin{equation}
\label{eq:ex_visco_11}
  \hat\sigma = E_0\hat\varepsilon + E_1\big(s^\alpha\hat\varepsilon-\fracint^{1-\alpha}\varepsilon(\nullng)\big)
\end{equation}
since $\alpha\in(0,1)$ and, therefore, $\nu=1$. Again, \eqref{eq:ex_visco_2} results from applying $\dds$ to \eqref{eq:ex_visco_11} and is rewritten here for completeness:
\begin{equation}
\label{eq:ex_visco_12}
  \hat\sigma^\prime = E_0\hat\varepsilon^\prime + E_1(\hat\varepsilon^\prime+\alpha s^{-1}\hat\varepsilon)s^\alpha.
\end{equation}
If the initial value $\fracint^{1-\alpha}\varepsilon(\nullng)$ is known, both equations can be combined such that $s^\alpha$ is eliminated. Similarly as in the case of homogeneous initial conditions, \eqref{eq:ex_visco_12} is obtained with $\hat\sigma+E_1\fracint^{1-\alpha}\varepsilon(\nullng)$ instead of $\hat\sigma$:
\begin{equation}
\label{eq:ex_visco_13}
  (\hat\varepsilon^\prime+\alpha s^{-1}\hat\varepsilon)\big(\hat\sigma+E_1\fracint^{1-\alpha}\varepsilon(\nullng)-E_0\hat\varepsilon\big) = \hat\varepsilon(\hat\sigma^\prime-E_0\hat\varepsilon^\prime).
\end{equation}
However, since the initial value $\fracint^{1-\alpha}\varepsilon(\nullng)$ is not simply the value of $\varepsilon$ at $t=0$ but is more complicated to interpret, it is only plausible to consider it as unknown. As discussed in the sequel, the unknown initial values (here only one) can either be eliminated by some algebraic manipulations or they can be treated as additional parameters and included in the identification. For the present example, the latter approach simply means that there are three parameters\footnote{The remaining parameter $E_1$ is again calculated from \eqref{eq:ex_visco_7}.} -- or products of parameters --  namely $\alpha$, $E_0$, and $E_1\fracint^{1-\alpha}\varepsilon(\nullng)$ to be identified in \eqref{eq:ex_visco_13}.

On the other hand, an elimination of an unknown initial value bears the advantage of being faced with less parameters to be identified. Clearly, \eqref{eq:ex_visco_12} does not depend on $\varepsilon(\nullng)$ anymore. So deriving once more w.r.t.\ $s$ allows for the elimination of $s^\alpha$:
\begin{multline}
  \big(\hat\varepsilon^\pprime+2\alpha s^{-1}\hat\varepsilon^\prime+\alpha(\alpha-1)s^{-2}\hat\varepsilon\big)(\hat\sigma^\prime-E_0\hat\varepsilon^\prime) \\
  = (\hat\varepsilon^\prime+\alpha s^{-1}\hat\varepsilon)(\hat\sigma^\pprime-E_0\hat\varepsilon^\pprime).
\end{multline}
The operational expressions in this equation can again be written as functions of time, which yields a nonlinear equation in $E_0$ and $\alpha$. In order to determine the unknown parameters either the nonlinear parameter problem can be solved or overparametrization can be used treating the coefficients $\alpha$, $\alpha^2$, $E_0$, $\alpha E_0$, and $\alpha^2 E_0$ as independent.

\subsubsection{Caputo fractional derivative}

The approach is a little different if Caputo's fractional derivative is used instead of the one due to Riemann and Liouville. In order to illustrate this, $\fracdiff^\alpha$ in \eqref{eq:ex_visco} is replaced by $\fracdiffc^\alpha$. The correspondence \eqref{eq:fracdiffc_op} then yields
\begin{equation}
\label{eq:ex_visco_21}
  \hat\sigma = E_0\hat\varepsilon + E_1\big(s^\alpha\hat\varepsilon-\varepsilon(\nullng)s^{\alpha-1}\big).
\end{equation}
As mentioned before, Caputo's definition bears the advantage of using initial values (here $\varepsilon(\nullng)$) that can easily be interpreted. While in this simple example the presumed knowledge of $\varepsilon(\nullng)$ does not impose any restriction, since $\varepsilon$ itself is considered as a known signal, for $\alpha>1$ initial values of derivatives come into play that are unknown in general. For that reason, and because in a real application measurements are noisy, $\varepsilon(0)$ is assumed to be unknown in the sequel.

Deriving \eqref{eq:ex_visco_21} w.r.t.\ $s$ gives
\begin{equation}
\label{eq:ex_visco_22}
  \hat\sigma^\prime = E_0\hat\varepsilon^\prime + E_1\big(s^\alpha(\hat\varepsilon^\prime+\alpha s^{-1}\hat\varepsilon)-(\alpha-1)s^{-1}\varepsilon(\nullng)s^{\alpha-1}\big).
\end{equation}
If the initial value is to be identified, eliminating $s^\alpha$ from \eqref{eq:ex_visco_22} using \eqref{eq:ex_visco_21},
\begin{multline}
  (\hat\sigma^\prime\hat\varepsilon-\hat\sigma\hat\varepsilon^\prime) = \big(\varepsilon(0)(1-\alpha)s^{-2}+\alpha s^{-1}\hat\varepsilon\big)(\hat\sigma-E_0\hat\varepsilon) \\
  + \varepsilon(0)s^{-1}(\hat\sigma^\prime-E_0\hat\varepsilon^\prime)
\end{multline}
is obtained which is an equation of the type \eqref{eq:hom_sys_4} and involves the parameters $\varepsilon(0)$, $E_0$, and $\alpha$.

Alternatively, the initial value $\varepsilon(0)$ can simply be eliminated by multiplying \eqref{eq:ex_visco_21} with $(1-\alpha)s^{-1}$ and adding it to \eqref{eq:ex_visco_22} multiplied with $s$:
\begin{equation}
  (1-\alpha)\hat\sigma + s\hat\sigma^\prime = E_0(1-\alpha)\hat\varepsilon+E_0s\hat\varepsilon^\prime+E_1s^\alpha(\hat\varepsilon+s\hat\varepsilon^\prime).
\end{equation}
By taking yet another derivative of this equation w.r.t.\ $s$ the factor $s^\alpha$ corresponding to the fractional derivative is replaced. Multiplying with $s^{-2}$, the resulting equation
\begin{multline}
  \hat\sigma^\prime\hat\varepsilon^\pprime-\hat\sigma^\pprime\hat\varepsilon^\prime + s^{-1}(\hat\sigma\hat\varepsilon^\pprime-\hat\sigma^\pprime\hat\varepsilon) + 2s^{-2}(\hat\sigma\hat\varepsilon^\prime-\hat\sigma^\prime\hat\varepsilon) \\
  = \alpha s^{-2}(\hat\sigma\hat\varepsilon^\prime-\hat\sigma^\prime\hat\varepsilon) - \alpha s^{-1}(2\hat\varepsilon^\prime+s^{-1}\hat\varepsilon)(\hat\sigma^\prime-E_0\hat\varepsilon^\prime) \\
  + \alpha s^{-1}\big(\hat\varepsilon^\pprime+\alpha s^{-1}\hat\varepsilon^\prime+(\alpha-1)s^{-2}\hat\varepsilon\big)(\hat\sigma-E_0\hat\varepsilon)
\end{multline}
is of the type \eqref{eq:hom_sys_4} and involves the parameters $E_0$ and $\alpha$.

\subsection{The general approach}

Based on the previous example and the discussions in section~\ref{sec:hom_gen} for the homogeneous case the general approach for an equation \eqref{eq:hom_sys} is rather obvious. That is why only a sketch is given of how to obtain relations of the form \eqref{eq:hom_sys_4}. For both the Riemann-Liouville and the Caputo definition the cases where initial conditions are eliminated and the one where they are identified are treated. Henceforth, it is assumed that an upper bound of all unknown (integer or fractional) orders $\bar\alpha_i$, $i=1,\dots,n+m$ is available\footnote{This assumption constitutes only a mild restriction since any (sufficiently large) finite number could be used as the upper bound.}.

\subsubsection{Riemann-Liouville fractional derivative}

For a known upper bound $\nu\in\Nset$ of the order of differentiation, the operational representation of \eqref{eq:hom_sys} (corresponding to \eqref{eq:hom_sys_2} for homogeneous initial conditions) can be written as
\begin{multline}
  \big(a_0(s)\hat y-b_0(s)\hat u\big)s^{\alpha_0} + \big(a_1(s)\hat y-b_1(s)\hat u\big)s^{\alpha_1} \\
  + \cdots + \big(a_r(s)\hat y-b_r(s)\hat u\big)s^{\alpha_r} = \textstyle\sum_{i=1}^\nu c_is^{i-1}
\end{multline}
where the coefficients $c_i$, $i=0,\dots,\nu$ depend on initial values of $u$ and $y$ and their derivatives as well as parameters $\bar a_i$, $0,\dots,n$ and $\bar b_i$, $i=0,\dots,m$. Applying $\dds^\nu$ to this equations eliminates all these initial values. Following the elimination of the fractional operators $s^{\alpha_i}$, $i=0,\dots,r$, a multiplication with $s^{-k}$ (with sufficiently large $k$) yields an equation of the form \eqref{eq:hom_sys_4} with $r+\nu$ instead of $r$.

On the other hand, identifying the initial values -- even though they lack physical meaning -- is almost identical to the procedure in section~\ref{sec:hom_gen}, the only difference being that \eqref{eq:hom_sys_3mat} is inhomogeneous with the right hand side of the equation being a vector of polynomials in $s$. Using the solution of this linear problem and substituting it in an additional equation involving the fractional orders $s^{\alpha_i}$ (generated from deriving once more w.r.t.\ $s$), again, a multiplication with a sufficiently large negative power of $s$ results in \eqref{eq:hom_sys_4}.

\subsubsection{Caputo fractional derivative}

The operational representation of \eqref{eq:hom_sys} with Caputo's derivation operator gives
\begin{multline}
\label{eq:inhom:cap}
  \big(a_0(s)\hat y-b_0(s)\hat u-\textstyle\sum_{i=1}^\nu c_{0,i}s^{-i}\big)s^{\alpha_0} \\
  + \cdots + \big(a_r(s)\hat y-b_r(s)\hat u-\textstyle\sum_{i=1}^\nu c_{r,i}s^{-i}\big)s^{\alpha_r} = 0
\end{multline}
for a known upper bound $\nu\in\Nset$. The coefficients $c_{i,j}$, $i=0,\dots,r$, $j=1,\dots,\nu$ can be cancelled by taking a sufficient number of derivatives w.r.t.\ $s$ and solving the resulting linear system of equations for these coefficients. Afterwards, an expression of the form \eqref{eq:hom_sys_4} can be obtained using $r$ further differentiations of this result and multiplying it with a sufficiently large negative power of $s$.

If, however, one is interested in the initial values, as in the homogeneous case, the fractional operators $s^{\alpha_i}$ in \eqref{eq:inhom:cap} can simply be eliminated by taking $r$ derivatives and calculating the solution of this linear system (cf.~\eqref{eq:hom_sys_3mat}).

\section{Linear partial differential equations with fractional derivatives}
\label{sec:pde}

% The solution in the case where there is a fractional derivative w.r.t.\ $z$ yields the Mittag-Leffler function. So far \emphng{I} have not found a helpful differential equation for this function (in order to eliminate it and go on with the identification).

The identification of fractional orders and model parameters discussed for lumped fractional models can be extended to fractional distributed parameter systems. For the integer order case, \cite{RudWoi08corsica} and \cite{geh12sysid} demonstrated that based on the knowledge of the solution of the corresponding operational ordinary differential equation a relation of the form \eqref{eq:hom_sys_4} between two boundary measurements can be obtained.

Here, the basic ideas of the identification method are demonstrated on an example. Since homogeneous initial conditions are assumed, without loss of generality, only the Riemann-Liouville fractional derivative is used.

The example considered is the fractional diffusion-wave equation (see e.g.~\cite{mainardi97frac,podlubny99book})
\begin{equation}
\label{eq:pde_ex}
  v^2\fracdiffz^2 u(z,t) = \fracdiff^\alpha u(z,t), \qquad (z,t)\in\Rset_+^2
\end{equation}
with $v\in\Rset$, where the derivative w.r.t.\ time is of fractional order $\alpha\in(0,2]$.
The factor $v^2$ is the diffusion coefficient for $\alpha=1$, and it is the squared wave propagation speed for $\alpha=2$. In the context of signal propagation as addressed in \cite{mainardi97frac} the boundary conditions are
\begin{equation}
\label{eq:pde_ex_bound}
  u(\nullng,t) = h(t), \qquad u(\infty,t) = 0.
\end{equation}
Additionally, as mentioned before, homogeneous initial conditions are assumed\footnote{Note that the number of initial conditions depends on the (unknown) fractional order $\alpha$. However, this fact is of no interest in the present discussion.}.

The operational representation of \eqref{eq:pde_ex} reads
\begin{equation}
  v^2\fracdiffz^2\hat u(z) = s^\alpha\hat u(z)
\end{equation}
and is a second order ordinary differential equation w.r.t.\ $z$ that can easily be solved. From the general solution
\begin{equation}
\label{eq:pde_ex_sol}
  \hat u(z) = \hat C_1\rme^{\frac{z}{v}s^{\alpha/2}} + \hat C_2\rme^{-\frac{z}{v}s^{\alpha/2}}
\end{equation}
using the (operational form) of the boundary conditions \eqref{eq:pde_ex_bound} yields operators $\hat C_1=0$ and $\hat C_2=\hat h$.

In order to identify the fractional order $\alpha$, apart from the knowledge of $h$, another measurement is required. If $u$ is measured at some point $z=L\neq0$, such that $g(t)=u(L,t)$ is known, \eqref{eq:pde_ex_sol} gives
\begin{equation}
\label{eq:pde_ex_1}
  \hat g = \hat h\rme^{-\frac{L}{v}s^{\alpha/2}}.
\end{equation}
Since the exponential function satisfies a first order differential equation it can be eliminated by deriving \eqref{eq:pde_ex_1},
\begin{equation}
\label{eq:pde_ex_2}
  \hat g^\prime = \left(\hat h^\prime-\tfrac{\alpha}{2}\tfrac{L}{v}s^{-1}s^{\alpha/2}\hat h\right)\rme^{-\frac{L}{v}s^{\alpha/2}}
\end{equation}
and replacing the exponential function using \eqref{eq:pde_ex_1}:
\begin{equation}
\label{eq:pde_ex_3}
  \hat g\hat h^\prime - \hat g^\prime\hat h = \tfrac{\alpha}{2}\tfrac{L}{v}s^{-1}s^{\alpha/2}\hat g\hat h.
\end{equation}
Finally, $s^{\alpha/2}$ is eliminated using the derivative of \eqref{eq:pde_ex_3} w.r.t.\ $s$. The interpretation of
\begin{multline}
\label{eq:pde_ex_4}
  \hat g^2\big((\hat h^\prime)^2-s^{-1}\hat h\hat h^\prime-\hat h\hat h^\pprime\big) - \hat h^2\big((\hat g^\prime)^2-s^{-1}\hat g\hat g^\prime-\hat g\hat g^\pprime\big) \\
  = \frac{\alpha}{2}s^{-1}\hat g\hat h(\hat g^\prime\hat h-\hat g\hat h^\prime)
\end{multline}
in terms of functions of time is omitted for brevity. It involves nested convolutions of known signals only. However, as can clearly be seen, \eqref{eq:pde_ex_4} constitutes a linear equation in $\alpha$. Based on an estimated value for the fractional order, the quotient\footnote{From a physical perspective it seems plausible that either the distance $L$ or the "wave propagation speed" $v$ have to be known in order to calculate the other one (based solely on the knowledge of $u(z,t)$ at two points $z$).} $L/v$ can directly be calculated from \eqref{eq:pde_ex_3} (since $s^{\alpha/2-1}$ already corresponds to a fractional integral for $\alpha\in(0,2)$).

% alternative: Taylor expansion for exponential function $\Rightarrow$ finite dimensional approximation

\section{Conclusion}

Using the algebraic framework presented model parameters and fractional (as well as integer) orders can be identified in linear fractional models. The exact relations obtained for the unknown quantities involve only convolutions of known (input and output) signals. Based on positive results using simulation data, the authors intend to further validate the identification algorithms using experimental data. In this context numerical aspects related to online implementation and robustness w.r.t.\ measurement noise should be investigated.

The example of the fractional diffusion-wave equation demonstrates, that the method is also applicable to (at least certain) linear fractional distributed parameter systems as well as delay systems with fractional derivatives. Beyond that, an extension towards the identification of the structure of linear systems, i.e., ordinary differential equations, seems possible where, both, the order of the equations and their coefficients can be estimated.
% extension towards the identification of linear systems structure, i.e., (integer) orders $n$, $m$ of and coefficients $a_i$, $b_i$ in differential equations
% \begin{equation}
%   y^{(n)} + a_{n-1}y^{(n-1)} + \cdots + a_1\dot y + a_0y = b_mu^{(m)} + \cdots + b_1\dot u + b_0u
% \end{equation}

% \begin{ack}                               % Place acknowledgements
% Partially supported by the Roman Senate.  % here.
% \end{ack}

% \bibliographystyle{alpha}        % Include this if you use bibtex 
% \bibliography{alles}           % and a bib file to produce the 
                                 % bibliography (preferred). The
                                 % correct style is generated by
                                 % Elsevier at the time of printing.

% \appendix
% \section{A summary of Latin grammar}    % Each appendix must have a short title.
% \section{Some Latin vocabulary}         % Sections and subsections are supported  
                                        % in the appendices.
\end{document}